\documentclass[12pt]{article}

\usepackage{graphicx}

\begin{document}

\title{$F$-split Galois representations are potentially abelian}

\author{Chandrashekhar Khare}

\date{}

\maketitle

\newtheorem{theorem}{Theorem}
\newtheorem{lemma}{Lemma}
\newtheorem{prop}{Proposition}
\newtheorem{cor}{Corollary}
\newtheorem{conj}{Conjecture}
\newtheorem{guess}{Guess}
\newtheorem{remark}{Remark}
\newtheorem{example}{Example}
\newtheorem{conjecture}{Conjecture}
\newtheorem{definition}{Definition}
\newtheorem{quest}{Question}
\newtheorem{ack}{Acknowledgemets}
\newcommand{\rhobar}{\overline{\rho}}
\newcommand{\Sha}{{\rm III}}

\noindent{\bf Abstract:} In this note we relate
the property of a semisimple $\ell$-adic Galois representation being
``$F$-split'' to its having abelian image.

\vspace{2mm}


\vspace{4mm}

If $E$ is an elliptic curve defined over a number field
$K$ then it gives rise to a {\it strictly compatible system}
$(\rho_{\ell})$ of 2-dimensional 
$\ell$-adic representations (see I-11 of [Se]) of the absolute Galois
group $G_K$ of $K$. This arises from the action of
$G_K$ on the $\ell$-adic Tate module ${\rm Ta}_{\ell}(E)$ of
$E$. 

This strictly compatible system has very different properties
when $E$ has CM and when $E$ does not have CM. For instance in the latter
case the image of $\rho_{\ell}$ is all of 
$GL_2({\bf Z}_{\ell})$ for almost all $\ell$ while in the former
case $(\rho_{\ell})$ is potentially abelian,
i.e., $\rho_{\ell}|_{G_L}$ is abelian for all $\ell$ for a 
fixed finite extension $L$ of $K$. These cases are also markedly different
from the point of view of the the distribution of the eigenvalues
of the Frobenii (see I-25 of [Se]).

Yet another notable difference is in the behaviour of
the compositum of the splitting fields 
of the {\it Frobenius polynomials} attached to $E$, i.e., characteristic polynomials 
of $\rho_{\ell}|_{G_L}({\rm Frob}_r)$ ($r$ prime to $\ell$) for
primes $r$ of $L$ at which $E_{/L}$ has good reduction: these
characteristic polynomials depend only
on $r$ and not on $\ell$, and here again $L$ is a fixed number field which
contains $K$. When $E$ does not have CM, 
this compositum is always an infinite extension of ${\bf Q}$ (see exercise
on page IV-13 of [Se]). If $E$ has CM this compositum is a finite
extension, whenever all the
endomorphisms (i.e., CM) of $E$ are defined over $L$.

This difference is the most pertinent to this short note.
We study here in a more abstract setting the relationship
between a compatible system of $\ell$-adic Galois representations being
potentially abelian and the nature of the field generated by the splitting fields
of the Frobenius polynomials.

\begin{definition} Let $E,F$ be number fields and fix embeddings of $E,F$ in each completion
of $\overline{\bf Q}$.
\begin{enumerate}
\item Consider $(\rho_{\lambda})$ a strictly compatible system
of $E$-rational, continuous, semisimple, $n$-dimensional, $\lambda$-adic representations $\rho_{\lambda}:G_K \rightarrow GL_n(E_{\lambda})$
for a number field $K$ and $\lambda$ running through the places of $E$ (see I-11 and I-13 of [Se]).
We say that $(\rho_{\lambda})$ is $F$-split  
if for almost all places $r$ of $L$ the characteristic
polynomial of $\rho_{\lambda}({\rm Frob}_r)$ (which is defined when 
$r$ and $\lambda$ are of coprime residue characteristics, and then
is independent of $\lambda$) 
{\it splits} over $F$. We say that $(\rho_{\lambda})$ is potentially $F$-split if
$(\rho_{\lambda}|_{G_L})$ is $F$-split for a finite extension $L$ of $K$. 
\item We say that a continuous, semisimple representation $\rho_{\lambda}:G_K \rightarrow GL_n(E_{\lambda})$
is $F$-split if the characteristic polynomial of $\rho_{\lambda}({\rm Frob}_r)$ of $L$ for all places $r$ of $L$ unramified
in $\rho_{\lambda}$ splits over $F$. 
We say that $\rho_{\lambda}$ is potentially $F$-split if
$\rho_{\lambda}|_{G_L}$ is $F$-split for a finite extension $L$ of $K$. (Note: By [KhRa], the set of
places $r$ of $L$ unramified in $\rho_{\lambda}$ has density 1.)
\end{enumerate}
\end{definition}

\begin{definition}
We say that either a strictly compatible system $(\rho_{\lambda})$ of $E$-rational, continuous, semisimple, $n$-dimensional 
$\lambda$-adic representations of $G_K$ or 
an $E$-rational, continuous, 
semisimple representation $\rho_{\lambda}:G_K \rightarrow GL_n(E_{\lambda})$ 
is abelian (resp., potentially abelian)
if each $\rho_{\lambda}$ has abelian image 
(resp., if for some finite extension $L$ of $K$
each $\rho_{\lambda}|_{G_L}$ has abelian image).
\end{definition}

\begin{theorem}
A strictly compatible system $(\rho_{\lambda})$ of $E$-rational, continuous, semisimple, $n$-dimensional 
$\lambda$-adic representations of $G_K$ is potentially $F$-split for some number field $F$ if and only if it 
is potentially abelian.
\end{theorem}

\noindent{\bf Proof:} We first prove the only if statement of the
theorem. By restriction of scalars, it is enough to prove that a
${\bf Q}$-rational, continuous, semisimple, $n$-dimensional 
$\ell$-adic system of representations of $G_K$ that is potentially 
$F$-split for some number field $F$ is potentially
abelian. (We reduce to this case for ease of comparison with [LP].)
Consider a prime $\ell_0$ that splits completely in $F$. 
Let $L$ be a finite extension of $K$ such that the compatible system
$(\rho_{\ell}|_{G_L})$ is $F$-split and such 
that the Zariski closure of the image $\rho_{\ell_0}(G_L)$ is
connected (such an 
$L$ exists!).  From the assumption of $F$-splitness
and the Cebotarev density 
theorem we deduce that the subgroup $\rho_{\ell_0}(G_L)$ of $GL_n({\bf Q}_{\ell_0})$ 
contains no non-split torus. Since $\rho_{\ell}$ is semisimple
this implies that the $\ell_0$-adic representation has toral image
(and 
in particular is abelian). This
together with the ${\bf Q}$-rationality of $\rho_{\ell_0}$,  and a consequence
of a result of Waldschmidt in transcendental number theory (see
Theorem 2 of [H]), 
implies that $\rho_{\ell_0}$ 
arises as the direct sum of 1-dimensional representations arising from 
algebraic Hecke characters $\chi_i$ of $L$, $i=1,\cdots,n$.
The $\chi_i$'s give rise 
(see II of [H]) to a strictly compatible system of (continuous,
semisimple, 
$n$-dimensional) $\lambda$-adic
representations. Comparing this with $(\rho_{\ell}|_{G_L})$, we deduce that $(\rho_{\ell}|_{G_L})$
itself ``arises'' from the direct sum of the 
algebraic Hecke characters $\chi_i$ proving the proposition. (Note that without 
using [He] it follows from Proposition 6.14 of [LP] that there is a finite extension $L$ of $K$ such that
for a density 1 set of primes $\ell$, $\rho_{\ell}|_{G_L}$ has abelian image, and thus $(\rho_{\ell}|_{G_L})$
is an ``abelian system'' with $\ell$ running through a density 1 set of primes. But we do not get the stronger assertion that
$(\rho_{\ell}|_{G_L})$ itself is an abelian system.)

We now prove the other direction of the statement of the theorem.
Consider a  strictly compatible system $(\rho_{\lambda})$ 
that is potentially abelian and consider an extension $L$ of $K$ such that
$(\rho_{\lambda}|_{G_L})$ is abelian. Choose any place $\lambda$ of $E$.
By appealing to [H] again, and using that $\rho_{\lambda}$ is rational over $E$, semisimple
and abelian we deduce that $\rho_{\lambda}|_{G_L}$
arises as the direct sum of 1-dimensional representations arising from 
algebraic Hecke characters $\chi_i$ of $L$.
Then by standard proprties of algebraic Hecke characters (see II of [H]), we conclude that 
$\rho_{\lambda}|_{G_L}$ is $F$-split for some number field $F$, which in its turn implies that
$(\rho_{\lambda}|_{G_L})$ is an $F$-split strictly compatible system.

\vspace{3mm}

\noindent{\bf Remark:} While it is true that an $E$-rational abelian, semisimple
representation $\rho_{\lambda}$ (or compatible system $(\rho_{\lambda})$) is always $F$-split for some number field $F$ (as
follows from the proof), it is not true that an $F$-split representation
$\rho_{\lambda}$ (or compatible system $(\rho_{\lambda})$) is abelian
(consider the ``constant'', compatible systems arising from Artin representations).

\vspace{3mm}

The case of a single,
$F$-split, $\lambda$-adic representation of $G_K$ we cannot answer even in the 
case when the image is in $GL_2({\bf Q}_{\ell})$ (if the completions of $F$ contain all the
quadratic extensions of ${\bf Q}_{\ell}$ 
we are at a loss how to proceed). We end with a question.

\begin{quest}
Is a continuous, semisimple, $F$-split $\lambda$-adic representation $\rho:G_K
\rightarrow GL_n(E_{\lambda})$ potentially abelian?
\end{quest}

\noindent{\bf Acknowledgements:} I thank Richard Pink, Dipendra Prasad and J-P.~Serre for
helpful correspondence: each came up with different proofs of the theorem of this note 
(independently of each other and of the author, when responding to messages posing 
mainly the question above). The terminology ``$F$-split'' was suggested by Richard Pink and 
J-P.~Serre pointed out the example in his book. I also thank David
Rohrlich for helpful comments that improved the exposition.

\section*{References}

\noindent [H] Henniart, G., {\it Repr\'esentations $\ell$-adiques 
ab\'eliennes}, in S\'eminaire  de Th\'eorie des
Nombres, Progress in Math. 22 (1982), 107--126, Birkhauser.

\vspace{3mm}


\noindent [KhRa] Khare, C., Rajan, C.~S.~, {\it The density
of ramified primes in semisimple $p$-adic Galois representations},
International Mathematics Research Notices  no. 12 (2001), 601--607.

\vspace{3mm}

\noindent [LP] Larsen, M., Pink, R., {\it On $\ell$-independence 
of algebraic monodromy groups in compatible systems of
representations}, Invent. Math. 107 (1992), 603--636. 

\vspace{3mm}

\noindent [S] Serre, J-P., {\it Abelian $\ell$-adic representations
and elliptic curves}, Addison-Wesley, 1989.

\vspace{3mm}

\noindent {\it Addresses of the author}: Dept of Math, Univ of Utah, 155 S 1400 E,
Salt lake City, UT 84112. e-mail address: shekhar@math.utah.edu
 
\noindent School of Mathematics, 
TIFR, Homi Bhabha Road, Mumbai 400 005, INDIA. 
e-mail addresses: shekhar@math.tifr.res.in

\end{document}